\documentclass[titlepage,11pt]{article}
\usepackage{latexsym}
\usepackage{amsfonts}
\usepackage{amsmath}
\usepackage{tikz, fp}
\usepackage{textcomp}
\usetikzlibrary{graphs}
\usepackage{float}
\usetikzlibrary{decorations.markings}
\usetikzlibrary{decorations,decorations.pathmorphing}
\newcommand\blackslug{\hbox{\hskip 1pt \vrule width 4pt height 8pt depth 1.5pt
        \hskip 1pt}}
\newcommand\bbox{\hfill \quad \blackslug \bigbreak}
\def\dd{\hbox{-}}

\def\ll{,\ldots,}
\newcommand{\R}{\mathbb{R}^+}

\title{Holes with hats and Erd\H{o}s-Hajnal}
\author{Maria Chudnovsky\thanks{This material is based upon work supported in part by the U. S. Army
Research Office under grant   number W911NF-16-1-0404, and supported by  NSF grant DMS 1763817.}\\
Princeton University, Princeton, NJ 08544
\\
\\
Paul Seymour\thanks{Supported by AFOSR grant
A9550-19-1-0187, and by NSF grant  DMS-1800053.}\\
Princeton University, Princeton, NJ 08544}

\date{May 1, 2020; revised \today}

\newtheorem{thm}{}[section]

\newcommand{\Proof}{\noindent{\bf Proof.}\ \ }
\newcommand{\vare}{\varepsilon}

\begin{document}
\maketitle
\begin{abstract}
A ``hole-with-hat'' in a graph $G$ is an induced subgraph of $G$ that consists of a cycle of length at least four, together 
with one further vertex that has exactly two neighbours in the cycle, adjacent to each other, and the ``house'' is the smallest, on
five vertices.
It is not known whether there exists $\vare>0$ such that every graph $G$ containing no house has a clique or stable set of 
cardinality at least $|G|^\vare$; this is one of the three smallest open cases of the Erd\H{o}s-Hajnal
conjecture and has been the subject of much study.

We prove that there exists $\vare>0$
such that every graph $G$ with no
hole-with-hat has a clique or stable set of cardinality at least $|G|^\vare$. 
\end{abstract}

\section{Introduction}
Graphs in this paper are finite and simple, and $|G|$ denotes the number of vertices of a graph $G$. 
A graph is {\em $H$-free} if it has no induced subgraph isomorphic to $H$.
The Erd\H{o}s-Hajnal conjecture~\cite{EH0,EH} asserts:
\begin{thm}\label{EHconj}
{\bf Conjecture:} For every graph $H$, there exists $\vare>0$ such that every $H$-free graph $G$ has a clique or stable set of
cardinality at least $|G|^\vare$.
\end{thm}
This has not yet been proved when $H$ is the five-vertex path $P_5$, and that problem motivates the work of this paper.
The complement of $P_5$ is the {\em house}, the graph consisting of a cycle of length four, together with one extra vertex
with two neighbours in the cycle, adjacent. By taking complements, we see that proving \ref{EHconj} when $H$ is the house is 
the same problem as proving it when $H=P_5$. The house is the smallest example of a ``hole-with-hat''.

\begin{figure}[H]
\centering

\begin{tikzpicture}[scale=0.8,auto=left]
\tikzstyle{every node}=[inner sep=1.5pt, fill=black,circle,draw]

\node (v0) at (0,2) {};
\node (v1) at (-1,0) {};
\node (v2) at (1,0) {};
\node (v3) at (-1,-2) {};
\node (v4) at (1,-2) {};

\draw (v0) -- (v1);
\draw (v0) -- (v2);
\draw (v1) -- (v2);
\draw (v2) -- (v4);
\draw (v3) -- (v4);
\draw (v1) -- (v3);

\end{tikzpicture}

\caption{A house.} \label{fig:house}
\end{figure}

A {\em hole} in a graph $G$ is an induced cycle of length at least four.  If $C$ is a hole in $G$, a
vertex $v\in V(G)\setminus V(C)$
is said to be a {\em hat} for $C$ if  $v$ has exactly two neighbours $x,y\in V(C)$, and $x,y$ are adjacent.
The subgraph induced on $V(C)\cup \{v\}$ is then said to be a {\em hole-with-hat} in $G$; and we say $G$ is {\em hole-with-hat-free}
if there is no hole-with-hat in $G$.

The main result of this paper is:
\begin{thm}\label{mainthm}
There exists $\vare>0$ such that for every hole-with-hat-free graph $G$, there is a clique or stable set in $G$
of cardinality at least $|G|^\vare$.
\end{thm}

Here are some earlier theorems with a similar nature:
\begin{itemize}
\item If $G$ contains no hole, there is a clique or stable set in $G$
of cardinality at least $|G|^{1/2}$. (This is immediate because such graphs are perfect.)
\item If $G$ contains no house and no hole of odd length, then again there is a clique or stable set in $G$
of cardinality at least $|G|^{1/2}$. (Again, because such graphs are perfect, a consequence of the ``strong perfect graph theorem''~\cite{SPGT}.)
\item For each $\ell>0$, there exists $\vare>0$ such that if $G$ contains no house and no hole of length at least $\ell$, then
$G$ has a clique or stable set of cardinality at least $|G|^\vare$. (This is a combination of a theorem of Bousquet, Lagoutte, and Thomass\'e~\cite{lagoutte} and a theorem of Bonamy, Bousquet and Thomass\'e~\cite{bonamy}.)
\end{itemize}

Working in structural graph theory, one always hopes to find a collection
of
non-crossing decompositions that together break the graph into simpler
pieces. That is because the existence of such a collection leads into the
well-understood area of tree-decompositions. However, such collections do
not often appear in the context of forbidden induced subgraphs.
Here we give a weakening  of this notion, that is still almost
as useful, and does work more often for induced subgraphs. 

Hole-with-hat-free graphs
typically admit a certain kind of separation, that we call a ``fracture''. (This is related to the ``amalgam'' decomposition of hole-with-hat-free graphs due to
Conforti,  Cornu\'ejols,  Kapoor and Vu\v{s}kovi\'c~\cite{conforti}.) A fracture is a certain kind of partition of the vertex set of our graph $G$ into three  parts (actually four parts, but we merge
two of them for this sketch) $A,B,C$, where $A,B$ are ``anticomplete'', that is, there are no edges between them.  
We call $A$ and $B$ the ``small'' and ``big'' sides of the fracture. (There is no symmetry
between $A$ and $B$ in the full definition of a fracture.) Let $S$ be the union of all small sides of fractures.
The graph $R$ obtained from $G$ by deleting $S$ does not admit a fracture with nonempty small side (because such a fracture
would extend to one in the whole graph $G$, and we would have deleted all its small side, including the part in $R$);
so $R$ has a very restricted type. We can base a proof on this, provided we can show that $R$ still contains a substantial 
part of $G$: in other words, that $S$ is not too big. And we could show this, if we could prove that:

\medskip

{\em  Every component of $S$
is anticomplete to the big side of some fracture of $G$. }

\medskip

This is where ``non-crossing'' would be useful.
If it were true that fractures form a set of non-crossing separations, then every component of $S$ would be a component of one
small side, and therefore anticomplete to the corresponding big side. This is not true, but we have a substitute: 
we can show that for any two fractures, if some component of the union of their small sides is not a component of 
either small side, then the two big sides are equal. It follows from this
that every component of $S$ is anticomplete to the big side of some fracture, which is what we needed.
A similar idea works in several other situations, and we hope to find further uses
for it in the future.

A graph $P$ is {\em perfect} if chromatic number equals clique number for every induced subgraph of $P$.
We denote the set of nonnegative real numbers by $\R$.
Let $G$ be a graph and $f:V(G)\rightarrow \R$; if $X\subseteq V(G)$, we define $f(X)=\sum_{v\in X}f(v)$, and if $P$ is an induced subgraph of $G$ we define 
$f(P)=f(V(P))$.
We say that $f$ is {\em good on $G$} if 
$f(P)\le 1$ for every perfect induced subgraph $P$ of $G$. Now let $\alpha\ge 1$. We denote by $f^\alpha$ the function $g$ on $G$
defined by $g(v)=(f(v))^{\alpha}$ for each $v\in V(G)$. Let us say that $G$ is {\em $\alpha$-narrow}
if $f^\alpha(G)\le 1$ for every good function $f$ on $G$. Here is a result of Chudnovsky and Safra~\cite{bull}, with a short proof by
Chudnovsky and Zwols~\cite{zwols}:

\begin{thm}\label{narrowtoEH}
If $G$ is $\alpha$-narrow then $G$ has a clique or stable set of cardinality at least $|G|^\vare$, where $\vare=1/(2\alpha)$.
\end{thm}
\Proof Let $P$ be a perfect induced subgraph of $G$ with as many vertices as possible, and let $p=|P|$. 
Let $f(v)=1/p$ for all $v\in V(G)$. Then $f$ is good on $G$, and so $f^\alpha(G)\le 1$, that is,
$p^{-\alpha}|G|\le 1 $, and so $p\ge |G|^{1/\alpha}=|G|^{2\vare}$. 
But $P$ is perfect, and so $P$, and hence $G$, has a clique or stable set of
cardinality at least $p^{1/2}\ge |G|^\vare$. This proves \ref{narrowtoEH}.~\bbox

In view of \ref{narrowtoEH}, the following implies \ref{mainthm}:
\begin{thm}\label{mainthm2}
There exists $\alpha\ge 1$ such that every hole-with-hat-free graph is $\alpha$-narrow.
\end{thm}

Let us say a pair $(G,f)$ is {\em $\alpha$-critical} if
$G$ is a graph and $f:V(G)\rightarrow \R$ is a good function on $G$, such that
\begin{itemize}
\item every proper induced subgraph of $G$ is $\alpha$-narrow; and
\item $f^\alpha(G)> 1$.
\end{itemize}

To prove \ref{mainthm2}, choose an appropriately large value of $\alpha$, and suppose (for a contradiction)
that \ref{mainthm2} is 
not satisfied, and look at a counterexample $G$ with as few vertices as possible. Hence there is a good function $f$ on $G$
with $f^\alpha(G)> 1$, and so $(G,f)$ is $\alpha$-critical. Consequently, \ref{mainthm2} can be reformulated as:
\begin{thm}\label{mainthm3}
There exists $\alpha\ge 1$ such that for every $\alpha$-critical pair $(G,f)$, there is a hole-with-hat in $G$.
\end{thm}

We will prove this at the end of the final section,
but let us sketch the proof now. We are proving $\alpha$-narrowness instead of proving the statement of
\ref{mainthm} directly, in order to handle homogeneous sets; so vertices will have non-negative weights, but
for this sketch the reader could assume that all
vertex-weights are one. If there are two disjoint anticomplete sets of vertices, that both have linear total
weight, then we win by induction; so we assume there are no two such sets. By a theorem of R\"odl, there is a subgraph
$X$ containing a linear fraction of the total weight of $G$, such that either $X$ is sparse (in a weighted sense) or its
complement is: and the second is impossible, by a theorem of Bonamy,  Bousquet and Thomass\'e. So the first holds. We look
at fractures of $X$. If $A,B,C$ is such a fracture, then $C$ has very small total weight, and so at least one of $A,B$ has big weight. But not both $A,B$ have linear
total weight, since they are anticomplete; and with some sleight of hand we can arrange that it is always the small
side $A$ that has small weight, and therefore that most of the weight of $X$ resides in $B$.

Let $S$ be the union of all
the small sides of fractures of $X$. By the remarkable fact that we described earlier, every component of $S$ is anticomplete
to some big side, and therefore has small weight; and so $S$ itself has small weight (because otherwise we could group
its components into two sets both with big weight). That means that deleting $S$ from $X$ gives a graph $R$ that still has
big weight. But every fracture in $R$ extends to a fracture in $X$ (this is another useful feature of fractures, and the
reason for using ``forcers'', which we do not explain here), and therefore $R$ has no fracture with
nonempty small side. Hence $R$ has a very restricted type, and in particular it is
$\alpha'$-narrow where $\alpha'$ is much less than $\alpha$; and it follows that $G$ itself is $\alpha$-narrow, which is what
we wanted to show. This completes the sketch.

\section{Complete pairs of sets}

Hole-with-hat-free graphs have some convenient structural properties, that we will prove next.
Let $A,B\subseteq V(G)$ be disjoint; we say they are {\em complete} to each other
if every vertex in $A$ is adjacent to every vertex in $B$, and
{\em anticomplete} if there are no edges between $A,B$.
We are concerned in this section with how the remainder of a hole-with-hat free graph can attach to a pair of sets of vertices 
that are complete to each other.
A graph is {\em anticonnected} if its complement is connected; and its {\em anticomponents}
are the complements of the components of its complement. If $X\subseteq V(G)$, we say that $X$ is {\em connected} if $G[X]$ is
connected, and {\em anticonnected} if $G[X]$ is anticonnected.
If $C\subseteq V(G)$, a vertex $v\in V(G)\setminus C$ is {\em mixed} on $C$ if $v$ is
neither complete not anticomplete to $C$.
We begin with:

\begin{thm}\label{wiggly1}
Let $G$ be a hole-with-hat-free graph. Let $C,D$ be disjoint anticonnected subsets of $V(G)$, complete to each other. Then no vertex of $V(G)\setminus (C\cup D)$ is both mixed on $C$ and mixed on $D$.
\end{thm}
\Proof
Suppose that $v\in V(G)\setminus (C\cup D)$ is both mixed on $C$ and mixed on $D$. Since $C$ is anticonnected, there exist nonadjacent
$c_1,c_2\in C$ such that $v$ is adjacent to $c_1$ and not to $c_2$; and choose $d_1,d_2\in D$ similarly. Then the subgraph
induced on $\{c_1,c_2,d_1,d_2,v\}$ is a house, contradicting that $G$ is hole-with-hat-free. This proves \ref{wiggly1}.~\bbox

\begin{thm}\label{wiggly2}
Let $G$ be a hole-with-hat-free graph. Let $C,D$ be disjoint subsets of $V(G)$, complete to each other,
such that $C$ is connected. Let $P$ be a connected subgraph of $G\setminus (C\cup D)$, such that
some vertex of $P$ has a neighbour in $C$, and
no vertex of $P$ is complete to $C$. For every $v\in V(P)$,
there exists $u\in V(P)$, mixed on $C$, such that every vertex in $D$ adjacent to $v$ is also adjacent to $u$.
\end{thm}
\Proof Suppose the claim is false, and choose a counterexample with $P$ minimal. Choose $v\in V(P)$ such that no vertex of
$P$ has a neighbour in $C$
and is adjacent to all neighbours of $v$ in $D$.
Choose $u\in V(P)$ mixed on $C$. By the minimality of $P$, $P$ is an induced path
between $u,v$, and $u$ is the only vertex of $P$ with a neighbour in $C$. Choose  $d\in D$ adjacent to $v$
and not to $u$. By the minimality of $P$, $v$ is the only vertex of $P$ adjacent to $d$. Since $C$ is connected, there exist
adjacent $c_1,c_2\in C$  such that $u$ is adjacent to $c_1$ and not to $c_2$. But then the subgraph induced on $V(P)\cup \{c_1,c_2,d\}$
is a hole-with-hat, a contradiction. This proves \ref{wiggly2}.~\bbox

\begin{thm}\label{wiggly3}
Let $G$ be a hole-with-hat-free graph. Let $C,D$ be disjoint subsets of $V(G)$, complete to each other,
such that $C$ is connected and $D$ is anticonnected. Let $P$ be a connected subgraph of $G\setminus (C\cup D)$, such that
some vertex of $P$ has a neighbour in $C$, and
no vertex of $P$ is complete to $C$. If some vertex of $P$ has a neighbour in $D$, then some vertex of $P$ is
mixed on $C$ and complete to $D$.
\end{thm}
\Proof
If some vertex of $P$ has a neighbour in $D$, then by \ref{wiggly2}, some vertex $v\in C$ is mixed on $C$ and has a neighbour in $D$, and
therefore by \ref{wiggly1}, $v$ is complete to $D$. This proves \ref{wiggly3}.~\bbox

\begin{thm}\label{wiggly4}
Let $G$ be a hole-with-hat-free graph. Let $C,D$ be disjoint nonempty subsets of $V(G)$, complete to each other, such that
$C$ is connected and anticonnected, and $D$ is anticonnected. Let $P$ be a connected subgraph of $G$ with $V(P)\cap (C\cup D)=\emptyset$,
such that some vertex of $P$ has a neighbour in $C$, and no vertex of $P$ is complete to $C$. Then no vertex of $P$ is mixed on $D$.
\end{thm}
\Proof
Suppose not, and choose $C,D,P$ not satisfying the theorem, with $P$ minimal.
Choose $u\in V(P)$ with a neighbour in $C$, and therefore mixed on $C$; and choose $v\in V(P)$ mixed on $D$.
From the minimality of $P$, it follows that $P$ is an induced path with ends $u,v$, and $u$ is the only vertex of $P$
with a neighbour in $C$, and no vertex of $P$ different from $v$ is mixed on $D$. By \ref{wiggly3}
$u$ is complete to $D$, and in particular $u\ne v$.
Let $u'$ be the neighbour  of $u$ in $P$. It follows that $C\cup \{u\}$ is connected and anticonnected,
and $u'$ is mixed on it; and this contradicts the minimality of $P$. This proves \ref{wiggly4}.~\bbox

\begin{thm}\label{wiggly5}
Let $G$ be a hole-with-hat-free graph. Let $C,D$ be disjoint nonempty subsets of $V(G)$, complete to each other, such that
$C$ is connected and anticonnected, and $D$ is connected and anticonnected.
Then there do not exist connected subgraphs
$P,Q$ of $G\setminus (C\cup D)$, with $V(P\cap Q)\ne \emptyset$, such that
\begin{itemize}
\item some vertex of $P$ has a neighbour in $C$, and no vertex of $P$ is complete to $C$;
\item some vertex of $Q$ has a neighbour in $D$, and no vertex of $Q$ is complete to $D$.
\end{itemize}
\end{thm}
\Proof
Suppose that such $P,Q$ exist, and choose $C,D,P,Q$ with $|V(P)|+|V(Q)|$ minimal. Choose $w\in V(P\cap Q)$, and choose $p\in V(P)$
with a neighbour in $C$, and $q\in V(Q)$ with a neighbour in $D$. From the minimality of $|V(P)|+|V(Q)|$, it follows
that $P$ is an induced path with ends $p,w$, and no vertex of $P$ different from $p$ has a neighbour in $C$; and similarly
for $Q$; and $V(P\cap Q)=\{w\}$. Suppose that $p$ is complete to $D$. Hence $p\notin V(Q)$, and so $p\ne w$.
Then $C\cup\{p\}$ is connected and anticonnected, and complete to $D$, and the two paths $P\setminus p$ and $Q$
contradict the minimality of $|V(P)|+|V(Q)|$. Thus $p$ is not complete to $D$. Since no other vertex of $P$ has a neighbour in $C$,
it follows from \ref{wiggly2} that no vertex of $P$ has a neighbour in $D$. Similarly no vertex of $Q$ has a neighbour in $C$.
In particular, no vertex of $P\cup Q$ is complete to $C$, contrary to \ref{wiggly4}. This proves \ref{wiggly5}.~\bbox

\section{Decomposing hole-with-hat-free graphs}

If $v\in V(G)$, we denote by $N(v)=N_G(v)$
the set of all neighbours of $v$ in $G$. If $N\subseteq V(G)$, we denote by $G[N]$
the induced subgraph with vertex set $N$.
A {\em weighted graph} is a pair $(G,w)$, where $G$ is a graph and $w:V(G)\rightarrow \R$ is a function, such that
$w(G)=1$. 
Let $\vare>0$. We say a weighted graph $(G,w)$ is {\em $\vare$-coherent} if
\begin{itemize}
\item for every $v\in V(G)$, $w(v)< \vare$;
\item for every $v\in V(G)$, $w(N_G(v))< \vare$; and
\item if $A,B\subseteq V(G)$ are disjoint and anticomplete then $\min (w(A), w(B))< \vare$.
\end{itemize}

First we need:
\begin{thm}\label{bigcomp}
Let $(G,w)$ be an $\vare$-coherent weighted graph. If $X\subseteq V(G)$ with $w(X)\ge 3\vare$, there is a component $Y$
of $G[X]$ with $w(Y)> w(X)-\vare$.
\end{thm}
\Proof Let $Z$ be a union of components of $G[X]$, minimal such that $w(Z)\ge \vare$. 
Since $X\setminus Z$ is anticomplete to $Z$, it follows that $w(X\setminus Z)<\vare$, and so $w(Z)> w(X)-\vare$. 
Choose a component $Y$ of $G\setminus X$ with $Y\subseteq Z$. 
From the minimality of $Z$, $w(Z\setminus Y)<\vare$, and so $w(Y)\ge ( w(X)-\vare)-\vare\ge \vare$, and therefore $Z=Y$
from the minimality of $Z$. But then $w(Y)=w(Z)> w(X)-\vare$. This proves \ref{bigcomp}.~\bbox

The component $Y$ of \ref{bigcomp} satisfies $w(Y)> w(X)-\vare\ge 2\vare$, and since the remainder of $G[X]$ is 
anticomplete to $Y$ and therefore has weight less than $\vare$, it follows that $Y$ is unique. We call $Y$ the {\em big component}
of $G[X]$.

If $A\subseteq V(G)$, each vertex in $V(G)\setminus A$ with a neighbour in $A$ is called an {\em attachment} of $A$.
Let $\vare>0$, with $5\vare\le 1$, and let $(G,w)$ be an $\vare$-coherent weighted graph. 
Let $C,D$ be disjoint subsets of $V(G)$, such that:
\begin{itemize}
\item $|C|\ge 2$, and $G[C]$ is connected and anticonnected;
\item $D\ne \emptyset$, and $D$ is the set of all vertices in $V(G)\setminus C$ that are complete to $C$; and
\item $C$ contains no attachment of the big component of $G\setminus (C\cup D)$.
\end{itemize}
(Note that since there  is a vertex in $D$ complete to $C$, it follows that $w(C)\le \vare$,
and similarly $w(D)\le \vare$.
Since $5\vare\le 1$, there is a big component
of $G\setminus (C\cup D)$.)
In these circumstances we call $(C,D)$ a {\em split} of $(G,w)$. As we shall see, splits are a useful
kind of decomposition in hole-with-hat-free graphs. 

Let us say a {\em forcer} is a graph $F$ with eight vertices $v_1\ll v_8$, where $v_1\dd v_2\dd v_3\dd v_4$
and $v_5\dd v_6\dd v_7\dd v_8$
are induced paths of $F$, and $\{v_1\ll v_4\}$ is complete to $\{v_5\ll v_8\}$. We call these two paths the {\em constituent paths}
of the forcer. A {\em forcer in $G$} means an induced subgraph
of $G$ that is a forcer, and $G$ is {\em forcer-free} if there is no forcer in $G$ .
Now we prove the main result of this section. It is a strengthening of the ``amalgam'' decomposition of hole-with-hat-free graphs due
to Conforti,  Cornu\'ejols,  Kapoor and Vu\v{s}kovi\'c~\cite{conforti}.

\begin{thm}\label{getsplit}
Let $\vare>0$, with $5\vare\le 1$, and let $(G,w)$ be a $\vare$-coherent weighted graph, where $G$
is hole-with-hat-free.
Let $F$ be a forcer in $G$. Then there is a split $(C,D)$ of $G$ such that $G[C], G[D]$ both contain a constituent path of $F$.
\end{thm}
\Proof
Let $F$ be a forcer, and let $P_1,P_2$ be the constituent paths of $F$.
Consequently there are disjoint subsets $X_1, X_2$ of $V(G)$, such that
\begin{itemize}
\item $V(P_1)\subseteq X_1$, and $X_1$ is connected and anticonnected;
\item $V(P_2)\subseteq X_2$, and $X_2$ is connected and anticonnected; and
\item $X_1, X_2$ are complete to one another.
\end{itemize}
Choose such $(X_1,X_2)$ maximal in the sense that there is no choice of $(X_1', X_2')$ satisfying the same conditions,
with $X_i\subseteq X_i'$ for $i = 1,2$ and $|X_1'\cup X_2'|>|X_1\cup X_2|$.
We call this property the {\em maximality} of $(X_1,X_2)$. Let $X_3$ be the set of all vertices in $V(G)\setminus (X_1\cup X_2)$
that are complete to $X_1\cup X_2$, and let $R=V(G)\setminus (X_1\cup X_2\cup X_3)$.
For $i=1,2$, let $R_i$ be the set of vertices in $R$ that are complete to $X_i$.
\\
\\
(1) {\em $R_1$ is anticomplete to $X_2$, and $R_2$ is anticomplete to $X_1$, and so $R_1\cap R_2=\emptyset$.}
\\
\\
Suppose that $v\in R_1$ has a neighbour in $X_2$, say. Since $v\notin X_3$, $v$ is mixed on $X_2$.
But then $X_2'=X_2\cup \{v\}$ is connected and anticonnected, and the pair $(X_1,X_2')$ violates the maximality of $(X_1,X_2)$.
This proves (1). 

\bigskip

For $i = 1,2$, let $S_i$ be the union of all components of $G[R\setminus (R_1\cup R_2)]$ that have an attachment in $X_i$.
Let $S_3=R\setminus (R_1\cup R_2\cup S_1\cup S_3)$.
\\
\\
(2) {\em $S_1\cap S_2=\emptyset$. Moreover, 
$S_1$ is anticomplete to $X_2\cup R_2\cup S_2$, and $S_2$ is anticomplete to $X_1\cup R_1\cup S_1$.}
\\
\\
By \ref{wiggly4}, $S_1\cap S_2=\emptyset$. By \ref{wiggly3}, $S_1$ is anticomplete to $R_2$, and $S_2$ is anticomplete to $R_1$.
This proves (2).

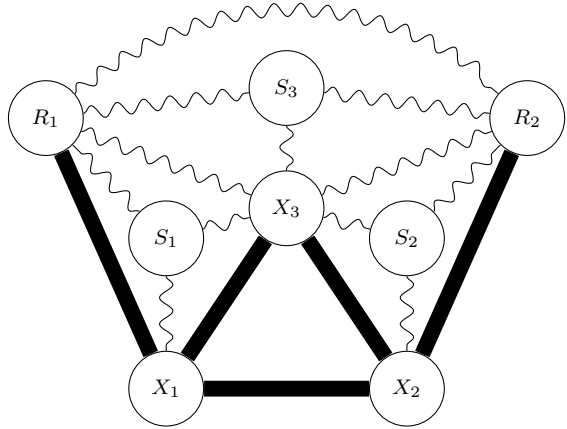
\begin{figure}[h!]
\centering

\begin{tikzpicture}[scale=0.8,auto=left]
\tikzstyle{every node}=[inner sep=10pt, fill=none,circle,draw]

\node (X3) at (0,1) {};
\node (S3) at (0,3) {};
\node (R1) at (-4,2.5) {};
\node (R2) at (4,2.5) {};
\node (X1) at (-2,-2) {};
\node (X2) at (2,-2) {};
\node (S1) at (-2,.5) {};
\node (S2) at (2,.5) {};

\draw[line width = 6pt] (R1) -- (X1);
\draw[line width = 6pt] (X1) -- (X2);
\draw[line width = 6pt] (X2) -- (R2);
\draw[line width = 6pt] (X2) -- (X3);
\draw[line width = 6pt] (X1) -- (X3);

\foreach \from/\to in {S1/R1, S1/X1,S1/X3,R1/X3,R1/S3,S3/X3,S3/R2,X3/R2,X3/S2,S2/R2,S2/X2}
\draw [decoration={snake}, decorate] (\from) -- (\to);

\draw [decoration={snake}, decorate] (R1) to [bend left =40](R2);

\tikzstyle{every node}=[]
\draw (R1) node []           {\scriptsize$R_1$};
\draw (R2) node []           {\scriptsize$R_2$};
\draw (X1) node []           {\scriptsize$X_1$};
\draw (S1) node []           {\scriptsize$S_1$};
\draw (X2) node []           {\scriptsize$X_2$};
\draw (S2) node []           {\scriptsize$S_2$};
\draw (S3) node []           {\scriptsize$S_3$};
\draw (X3) node []           {\scriptsize$X_3$};

\end{tikzpicture}

\caption{Thick lines indicate complete pairs, and wiggly lines indicate possible edges.} \label{fig:getsplit}
\end{figure}

Thus, in summary, the sets $X_1,X_2,X_3,R_1,R_2,S_1,S_2,S_3$
are pairwise disjoint and have union $V(G)$. The pairs
$$(X_1,X_2), (X_1,X_3), (X_2,X_3), (R_1,X_1), (R_2,X_2)$$
are complete to each other; the pairs
$$(R_1,X_2), (R_2,X_1), (S_1,X_2), (S_2,X_1), (S_3,X_1), (S_3,X_2), (S_1,R_2), 
(S_2,R_1), (S_1,S_2), (S_1,S_3),(S_2,S_3)$$
are anticomplete; and there may be edges between the pairs not listed.
Every component of $G[S_i]$ has an attachment in $X_i$ for $i = 1,2$.

Define $T=X_1\cup X_2\cup X_3\cup R_1\cup R_2$. Choose $x_1\in X_1$ and $x_2\in X_2$. Then every vertex in $T$ 
is adjacent to one of $x_1,x_2$,
and so $w(T)\le 2\vare$. Hence $G\setminus T$ has a big component $Y$, and since the sets $S_1,S_2,S_3$ are pairwise anticomplete,
we may assume that $Y$ is disjoint from $S_1$, by exchanging $X_1,X_2$ if necessary.
But then $(X_1,X_2\cup X_3\cup R_1)$ is a split of $G$ satisfying the theorem. This proves 
\ref{getsplit}.~\bbox

Let us say a split $(C,D)$ of $G$ is {\em optimal} if there is no split $(C',D')$
with $C\subseteq C'$ and $C'\ne C$. Let $(C,D)$ be an optimal split. Let $A$ be the union of all components of 
$G\setminus (C\cup D)$ that have an attachment in $C$; and let $B$ be the union of all other components of $G\setminus (C\cup D)$
(including the big component). Let us call $(A,C,D,B)$ a {\em fracture} of $G$.
(Note that there are no edges between $B$ and $A\cup C$ but there may well be edges between $A$ and $D$. Also, $B\ne \emptyset$,
since it contains the big component of $G\setminus (C\cup D)$, but $A$ might be empty.)
From \ref{getsplit} we have immediately:

\begin{thm}\label{getfracture}
Let $\vare>0$, with $5\vare\le 1$, and let $(G,w)$ be an $\vare$-coherent weighted graph, where
$G$ is hole-with-hat-free.
Let $F$ be a forcer in $G$. Then there is a fracture $(A,C,D,B)$ of $G$ such that $G[C]$ contains a  constituent path of $F$.
\end{thm}

We need some observations about fractures.
\begin{thm}\label{fracture}
Let $\vare>0$, with $5\vare\le 1$, and let $(G,w)$ be an $\vare$-coherent weighted graph, where
$G$ is hole-with-hat-free.
Let $(A,C,D,B)$ be a fracture of $G$. 
\begin{itemize}
\item For each $a\in A$, there is an attachment of the big component of $G\setminus (C\cup D)$ that is nonadjacent to~$a$.
\item For each $a\in A$, and every anticomponent $X$ of $G[D]$, $a$ is not mixed on~$X$.
\end{itemize}
\end{thm}
\Proof
Let $Y$ be the big component of $G\setminus (C\cup D)$; thus $Y\subseteq B$, and all its attachments belong to $D$.
Suppose that $a\in A$, and $a$ is adjacent to every vertex of $D$ that has a neighbour in $Y$.
Let $P$ be the component of $G[A]$ that contains $a$; then some attachment of $P$ belongs to $C$. By \ref{wiggly2},
we may choose $v\in P$ mixed on $C$, such that every vertex in $D$ adjacent to $a$ is also adjacent to $v$.
Let $D'$ be the set of all neighbours of $v$ in $D$. Then $(C\cup \{v\}, D')$ is a split 
(because $D'$ contains all attachments of $Y$), contradicting that $(C,D)$ is optimal. This proves the first assertion.

Now suppose that $a\in A$ is mixed on an anticomponent $X$ of $G[D]$. Let
$P$ be the component of $G[A]$ that contains $a$. Choose $v\in P$ mixed on $C$; then $P$ contradicts \ref{wiggly4} applied to the
connected anticonnected set $C$ and the anticonnected set $X$. This proves \ref{fracture}.~\bbox

\section{Multiple fractures}
A fracture $(A,C,D,B)$ of $G$ is a kind of separation of $G$, because deleting $C\cup D$ disconnects $A$ from $B$. 
(But $A$ might be empty.) Also the order of this separation is small, since $w(C\cup D)\le 2\vare$ in the usual notation.
 It would be nice if these separations did not ``cross'', so that they give us a tree-decomposition of $G$, but that is not true. 
Nevertheless, something like that is true, as we see in this section.

\begin{thm}\label{crossing}
Let $\vare>0$, with $6\vare\le  1$, and let $(G,w)$ be an $\vare$-coherent weighted graph, where
$G$ is hole-with-hat-free. Let $(A,C,D,B)$ and $(A', C', D', B')$ be fractures in $G$. Then either
\begin{itemize}
\item every connected subgraph of $G[A\cup A']$ is contained in one of $A$,  $A'$; or
\item the big component of $G\setminus (C\cup D)$ equals the big component of $G\setminus (C'\cup D')$.
\end{itemize}
\end{thm}
\Proof
We suggest that, to follow this argument, the reader imagine a $4\times 4$ matrix with rows labelled $A,C,D,B$ and columns 
$A', C', D', B'$. We remind the reader that $C$ is complete to $D$, and $A\cup C$ is anticomplete to $B$, and the same for $(A',C',D',B')$.

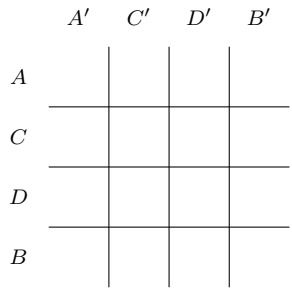
\begin{figure}[H]
\centering

\begin{tikzpicture}[scale=0.8,auto=left]
\tikzstyle{every node}=[]

\draw [-] (1,-2) -- (1,2);
\draw [-] (0,-2) -- (0,2);
\draw [-] (-1,-2) -- (-1,2);
\draw [-] (-2,1) -- (2,1);
\draw [-] (-2,0) -- (2,0);
\draw [-] (-2,-1) -- (2,-1);

\node at (-2.5,1.5) {\scriptsize $A$};
\node at (-2.5,0.5) {\scriptsize $C$};
\node at (-2.5,-0.5) {\scriptsize $D$};
\node at (-2.5,-1.5) {\scriptsize $B$};
\node at (-1.5,2.5) {\scriptsize $A'$};
\node at (-.5,2.5) {\scriptsize $C'$};
\node at (.5,2.5) {\scriptsize $D'$};
\node at (1.5,2.5) {\scriptsize $B'$};

\end{tikzpicture}

\caption{Two fractures.} \label{fig:matrix1}
\end{figure}

Let $Y$ be the big component of $G\setminus (C\cup D)$, and define $Y'$ similarly.
Since $w(Y), w(Y')> 1-3\vare\ge 1/2$, it follows that $Y\cap Y'\ne \emptyset$, and since $Y\subseteq B$ and $Y'\subseteq B'$, 
we deduce that $Y\cap Y'\cap B\cap B'\ne \emptyset$.
\\
\\
(1) {\em If $C\cap B'\ne \emptyset$ then $D\cap (A'\cup C')=\emptyset$, and if 
$B\cap C'\ne \emptyset$ then $(A\cup C)\cap D'=\emptyset$.}
\\
\\
Let $u\in C\cap B'$. If $v\in D\cap (A'\cup C')$, then $v$ is adjacent to $u$ (because $C$ is complete to $D$), and yet 
$v$ is nonadjacent to $u$ (because $A'\cup C', D'$ are anticomplete), a contradiction. This proves the first statement, and the second
follows by symmetry.
\\
\\
(2) {\em We may assume that $A\cap (C'\cup D')\ne \emptyset$, and $(C\cup D)\cap A'\ne \emptyset$, and at least one of 
$B\cap D', D\cap B'$ is nonempty.}
\\
\\
If $A\cap (C'\cup D')= \emptyset$, then the first outcome of the theorem holds, and similarly if $(C\cup D)\cap A'=\emptyset$.
If $B\cap D', D\cap B'$ are both empty, then $Y,Y'\subseteq B\cap B'$, and so $Y=Y'$ and the second outcome of the theorem holds.
This proves (2).

\bigskip
From the third assertion of (2) and symmetry, we may assume that $B\cap D'\ne \emptyset$. By (1), $(A\cup C)\cap C'=\emptyset$.
From (2), $A\cap D'\ne \emptyset$; so by (1) $B\cap C'=\emptyset$. Hence $D\cap C'\ne \emptyset$, because $C'\ne \emptyset$.
Every vertex in
$C\cap A'$ is complete to $D\cap C'$ and hence to $C'$; but no vertex in $A'$ is complete to $C'$ from the definition
of a fracture, and so $C\cap A'=\emptyset$. By (2), $D\cap A'\ne \emptyset$. By (1), $C\cap B'=\emptyset$, and so $C\cap D'\ne \emptyset$.

\begin{figure}[H]
\centering

\begin{tikzpicture}[scale=0.8,auto=left]
\tikzstyle{every node}=[]

\draw [-] (1,-2) -- (1,2);
\draw [-] (0,-2) -- (0,2);
\draw [-] (-1,-2) -- (-1,2);
\draw [-] (-2,1) -- (2,1);
\draw [-] (-2,0) -- (2,0);
\draw [-] (-2,-1) -- (2,-1);

\node at (-2.5,1.5) {\scriptsize $A$};
\node at (-2.5,0.5) {\scriptsize $C$};
\node at (-2.5,-0.5) {\scriptsize $D$};
\node at (-2.5,-1.5) {\scriptsize $B$};
\node at (-1.5,2.5) {\scriptsize $A'$};
\node at (-.5,2.5) {\scriptsize $C'$};
\node at (.5,2.5) {\scriptsize $D'$};
\node at (1.5,2.5) {\scriptsize $B'$};

\node at (-1.5,1.5) {\scriptsize $?$};
\node at (-1.5,.5) {\scriptsize $\emptyset$};
\node at (-1.5,-1.5) {\scriptsize $?$};
\node at (-.5,1.5) {\scriptsize $\emptyset$};
\node at (-.5,.5) {\scriptsize $\emptyset$};
\node at (-.5,-1.5) {\scriptsize $\emptyset$};
\node at (.5,-.5) {\scriptsize $?$};
\node at (1.5,1.5) {\scriptsize $?$};
\node at (1.5,.5) {\scriptsize $\emptyset$};
\node at (1.5,-.5) {\scriptsize $?$};

\tikzstyle{every node}=[inner sep=1.5pt, fill=black,circle,draw]
\node at (-1.5,-.5) {};
\node at (-.5,-.5) {};
\node at (.5,1.5) {};
\node at (.5,-1.5) {};
\node at (1.5,-1.5) {};
\node at (.5,.5) {};

\end{tikzpicture}
\caption{A solid dot means a nonempty set, and $?$ means we don't know.} \label{fig:matrix2}

\end{figure}

Since $C=C\cap D'$ is anticonnected, and every vertex
in $A$ has a nonneighbour in $C$, and every vertex in $B\cap D'$ has a nonneighbour in $A\cap D'$, it follows that
$(A\cup B\cup C)\cap D'$ is anticonnected. But each vertex in $D\cap A'$ has a neighbour in $(A\cup B\cup C)\cap D'$ (namely, in
$C\cap D'$), and by \ref{fracture}, it follows that $D\cap A'$ is complete to $(A\cup B\cup C)\cap D'$.
Similarly, since $D\cap (A'\cup B'\cup C')$ is anticonnected, it follows that $A\cap D'$ is complete to 
$D\cap (A'\cup B'\cup C')$.
(Thus we almost have symmetry between $(A,C,D,B)$ and $(A', C', D', B')$; but not quite, because we do not know that $D\cap B'\ne \emptyset$.)

Let $Q$ be the set of vertices in $D\cap D'$ that are not complete to $D\cap A'$, and let $Q'$ be the set of vertices in 
$D\cap D'$ that are not complete to $ A\cap D'$. Let $R=(D\cap D')\setminus (Q\cup Q')$.
\\
\\
(3) {\em $Q\cap Q'=\emptyset$, and $Q,Q', R$ are pairwise complete.}
\\
\\
Since there is a vertex in $A\cap D'$ and it is complete to $D\cap A'$, \ref{fracture} implies that $A\cap D'$ is complete to
$Q$; and so $Q\cap Q'=\emptyset$. If $u\in Q$ and $v\in Q'\cup R$, there is a vertex in $D\cap A'$ adjacent to $v$ and not to $u$;
so \ref{fracture} implies that $u,v$ are adjacent. Thus $Q$ is complete to $Q'\cup R$, and similarly $Q'$ is complete to $R$.
This proves (3).
\\
\\
(4) {\em $A'\cap B$ is anticomplete to $Q'\cup (B\cap D')$, and $B'\cap A$ is anticomplete to $Q\cup (D\cap B')$.}
\\
\\
Each vertex in $A'\cap B$ is anticomplete to $A\cap D'$, and so by \ref{fracture}, also anticomplete to $Q'\cup (B\cap D')$. Similarly
$B'\cap A$ is anticomplete to $Q\cup (D\cap B')$. This proves (4).

\bigskip

Choose $v\in D\cap A'$; then by \ref{fracture}, there is an
attachment $q$ of $Y'$ nonadjacent to $v$. Since $v$ is adjacent to all vertices of $D'\setminus Q$, it follows that
$q\in Q$.
Similarly there is an attachment $q'$ of $Y$ with $q'\in Q'$. Let $X=((A\cup C)\cap D')\cup Q'$,
and $X'=(D\cap (A'\cup C'))\cup Q$.
Then $X, X'$ are disjoint, and complete to each other, and each of them is both connected and anticonnected.
Now some vertex of $Y$ is adjacent to $q'$, and so has a neighbour in $X$; and no vertex of $Y$ is complete to $X$
(because $Y\subseteq B\cap (B'\cup D')$, since there are no edges between $B\cap A'$ and $B\cap D'$). Similarly
some vertex of $Y'$ has a neighbour in $X'$, and no vertex of $Y'$ is complete to $X'$. Since $Y\cap Y'$ is non-null, 
this contradicts \ref{wiggly5}.~\bbox

\begin{thm}\label{bigbag}
Let $\vare>0$, with $6\vare\le  1$, and let $(G,w)$ be an $\vare$-coherent weighted graph, where
$G$ is hole-with-hat-free. Let $\mathcal{F}$ be the set of all fractures of $G$, and let $\mathcal{A}$ be the union of
all the sets $A$ for $(A,C,D,B)\in \mathcal{F}$. Then $w(\mathcal{A})< 3\vare$.
\end{thm}
\Proof 
Let $Z$ be the vertex set of a component of $G[\mathcal{A}]$.
For each $(A,C,D,B)\in \mathcal{F}$, we call each component of $G[A]$ a {\em piece}; let $H$ be the set of all maximal
pieces (taken over
all $(A,C,D,B)\in \mathcal{F}$). Thus $Z$ can be expressed as the union of vertex sets of maximal pieces. For each maximal piece $X$,
let $(A,C,D,B)\in \mathcal{F}$ such that $X$ is a component of $G[A]$, and let $Y$ be the big component of 
$G\setminus (C\cup D)$; we call
$Y$ the {\em fulcrum} of $X$. (There may be more than one choice of $(A,C,D,B)\in \mathcal{F}$ for a given set $X$, and correspondingly
more than one choice of fulcrum: choose one, arbitrarily).

We observe:
\\
\\
(1) {\em If $X, X'$ are maximal pieces such that either $V(X\cap X')\ne \emptyset$, or $X$ is not anticomplete to $X'$, 
then $X,X'$ have the same fulcrum.}
\\
\\
Suppose not. Let $X$ be a component of $G[A]$ where $(A,C,D,B)\in \mathcal{F}$, and define $(A',C', D', B')$ similarly. 
By \ref{crossing}, it follows that every connected subgraph of $G[A\cup A']$ is a subgraph of one of $G[A], G[A']$, and in particular
the connected subgraph induced on $V(X)\cup V(X')$ is a subgraph of one 
of $G[A], G[A']$, say $G[A]$. But $X$ is a component of $G[A]$, so $V(X)=V(X\cup X')$, contradicting that $X'$ is a maximal piece.
This proves (1).

\bigskip

Choose a connected subgraph $H$ of $G[Z]$, maximal such that $V(H)$ is the union of maximal pieces all with the same fulcrum $Y$.
Suppose that $V(H)\ne Z$. Since $G[Z]$ is connected, there is a vertex $v_1\in Z\setminus V(H)$ with a neighbour $v_2\in V(H)$. Choose
a maximal piece $X_1$ containing $v_1$, and a maximal piece $X_2$ containing $v_2$ with fulcrum $Y$. By (1), $X_1$ has fulcrum $Y$,
contrary to the maximality of $H$. Thus $V(H)=Z$, and so $Y$ is anticomplete to $Z$. Since $w(Y)\ge \vare$, it follows that $w(Z)<\vare$.
Since this holds for each component of $G[\mathcal{A}]$, \ref{bigcomp} implies that $w(\mathcal{A})<3\vare$. This proves \ref{bigbag}.~\bbox

Let us say $X\subseteq V(G)$ is a {\em homogeneous set} of $G$ if for every vertex $v\in V(G)\setminus X$, either $v$ is complete
or anticomplete to $X$. Let $G$ be a graph; we say that $G$ is {\em guarded} if for every forcer $F$ in $G$,
there is a homogeneous set $X$ of $G$ with $X\ne V(G)$ such that $G[X]$ contains a constituent path of $F$.

\begin{thm}\label{flat}
Let $\vare>0$, with $6\vare\le  1$, and let $(G,w)$ be an $\vare$-coherent weighted graph, where
$G$ is hole-with-hat-free. Then there exists 
$Z\subseteq V(G)$ with $|Z|>1$, such that
$G[Z]$ is connected and guarded, and
$w(Z)>1-4\vare$.
\end{thm}
\Proof
Define $\mathcal{F}, \mathcal{A}$ as in \ref{bigbag}, and let $W=V(G)\setminus \mathcal{A}$. By \ref{bigbag}, 
$w(W)> 1-3\vare\ge 3\vare$. By \ref{bigcomp}, $G[W]$ has a big component, with vertex set $Z$ say, where $w(Z)\ge 1-4\vare$.
Hence $|Z|>1$, since $w(v)\le \vare< 1-4\vare$ for each vertex $v$.
Let $F$ be a forcer in $G[Z]$. Then by \ref{getfracture}, there is a fracture $(A,C,D,B)$ of $G$ such that $|V(F)\cap C|\ge 4$.
Let $X=C\cap Z$. Since $C$ is a homogeneous set of $G\setminus A$, it follows that $X$ is a homogeneous set of $G[Z]$, and it contains
a constituent path of $F$. This proves \ref{flat}.~\bbox

\section{$\alpha$-critical pairs}

In this section we explore the properties of $\alpha$-critical pairs, and combine these results with \ref{flat} to prove \ref{mainthm3}. 
\begin{thm}\label{smalldeg}
Let $\alpha\ge 2$, and let $(G,f)$ be $\alpha$-critical. Then $f(w)< 1-4^{-1/\alpha}$ for each $w\in V(G)$.
\end{thm}
\Proof
Let $w\in V(G)$, and let $c=f(w)$.
Let $N=N_G(w)$, and let $M=V(G)\setminus (N\cup \{w\})$. Since $(G,f)$ is $\alpha$-critical,
it follows that $G[N]$ is $\alpha$-narrow, and so is $G[M]$. Let $p$ be the maximum of $f(P)$ over all perfect induced subgraphs
of $G[N]$, and let $q$ be the maximum of $f(Q)$ over all perfect induced subgraphs $Q$
of $G[M]$. We claim that $f^\alpha (M)\le p^\alpha$. If $f(v)=0$ for every $v\in N$ then the statement is true, so we may assume
that $f(v)>0$ for some $v\in N$, and hence $p>0$.
So the function $f(v)/p\;(v\in N)$ is a good function on $G[N]$, and since $G[N]$ is $\alpha$-narrow, we deduce that
$f^\alpha (N)\le p^\alpha$. Similarly $f^\alpha(M)\le q^\alpha$.

But if $P$ is a perfect induced subgraph of $G[N]$ then $G[V(P)\cup \{w\}]$ is perfect, and therefore $f(V(P)\cup \{w\})\le 1$;
and so
$p\le 1-c$, and similarly $q\le 1-c$. Thus
$$1<f^\alpha(G)=f^\alpha(N)+f^\alpha(M)+f^\alpha(w)\le p^\alpha+q^\alpha+f^\alpha(w)\le 2(1-c)^\alpha+c^\alpha.$$

Now for $0\le x\le 1$, the function $g(x)=2(1-x)^\alpha+x^\alpha$ has the value $1$ when $x=1$, and its value increases with $x$
for $2/3\le x\le 1$, since $\alpha\ge 2$ (as can be seen by taking the derivative). Thus $g(x)\le 1$ for $2/3\le x\le 1$.
Since $g(c)>1$, it follows that $c<2/3$, and so $c^\alpha\le 1/2$; and consequently $2(1-c)^\alpha> 1/2$, that is,
$c<1-4^{-1/\alpha}$. This proves \ref{smalldeg}.~\bbox

\begin{thm}\label{strongEH}
Let $\alpha\ge 1$, and let $(G,f)$ be $\alpha$-critical. Let $A,B\subseteq V(G)$ be disjoint and either complete or anticomplete.
Then not both $f^\alpha(A), f^\alpha(B)> 2^{-\alpha}$.
\end{thm}
\Proof
Let $P$ be a perfect induced subgraph of $G[A]$ with $f(P)$ maximum, and choose $Q$ in $G[B]$ similarly. Since $P\cup Q$
is a perfect induced subgraph of $G$, it follows that $f(P)+f(Q)\le 1$, and from the symmetry we may assume that $f(P)\le 1/2$.
We may also assume that $f(P)>0$, $f(P)=p$ say, and so $f(v)/p\;(v\in A)$ is a good function on $G[A]$; and we may assume that
$Y\ne \emptyset$, and so $G[A]$ is $\alpha$-narrow, and consequently
$f^\alpha(A)\le p^\alpha\le 2^{-\alpha}$. This proves \ref{strongEH}.~\bbox

Next we need the following consequence of a theorem of R\"odl~\cite{rodl}:
\begin{thm}\label{rodl}
For all $\vare>0$ and every graph $H$, there exists $\delta > 0$ such that for every $H$-free graph
$G$, there is a subset $X\subseteq V(G)$ with $|X|\ge \delta |V(G)|$ such that
one of $G[X], \overline{G}[X]$ has maximum degree less than $\vare |X|$.
\end{thm}
We also need the following theorem of Bousquet, Lagoutte and Thomass\'e~\cite{lagoutte}:
\begin{thm}\label{lagoutte}
For every path $H$, there exists $\vare>0$ such that for every $H$-free graph $G$ with $|G|>1$, either some vertex of $G$ has degree
at least $\vare|G|$, or there are disjoint anticomplete subsets $A,B\subseteq V(G)$ with $|A|,|B|\ge \vare|G|$.
\end{thm}
We recall that the {\em house} is the complement of $P_5$. Let us say $G$ is {\em house-free} if $G$ contains no house.
\begin{thm}\label{nohouse}
For all $\vare>0$ there exists $\delta>0$ such that, if $G$ is house-free and $|G|>1$, then
either
\begin{itemize}
\item there are disjoint sets $A,B\subseteq V(G)$, complete to each other, with $|A|,|B|\ge \vare\delta |G|$, or
\item there exists $X\subseteq V(G)$ with $|X|\ge \delta|G|$ such that $G[X]$ has maximum degree less than $\vare|X|$.
\end{itemize}
\end{thm}
\Proof
Choose $\vare'>0$ such that \ref{lagoutte} holds with $H, \vare$ replaced by $P_5, \vare'$ respectively. Now let $\vare>0$;
we must show that there exists $\delta$ as in the theorem. Thus we may assume that $\vare\le \vare'$, by reducing $\vare$ if
necessary. Choose $\delta$ as in \ref{rodl}. Now let $G$ be a house-free graph with $|G|>1$. The complement
$\overline{G}$ of $G$ is $P_5$-free, and so by \ref{rodl}, there is a subset $X\subseteq V(G)$ with $|X|\ge \delta |V(G)|$ such that
one of $G[X], \overline{G}[X]$ has maximum degree less than $\vare |X|$.

If  $G[X]$ has maximum degree less than $\vare |X|$ then the theorem holds, so we assume that $G[X]$ has maximum degree at least 
$\vare |X|$ (and so $|X|>1$), and therefore
$\overline{G}[X]$ has maximum degree less than $\vare |X|$.
By \ref{lagoutte} applied to $\overline{G}[X]$,
there are disjoint anticomplete (in $\overline{G}$) subsets $A,B\subseteq X$ with $|A|,|B|\ge \vare|X|$. But then
$A$ is complete to $B$ in $G$, and $|A|,|B|\ge \vare|X|\ge \delta\vare |G|$. This proves \ref{nohouse}.~\bbox

From \ref{nohouse} we deduce:
\begin{thm}\label{fnohouse}
Let $\vare>0$, and choose $\delta>0$ satisfying \ref{nohouse}.
Let $\alpha\ge 1$, such that  $\vare\delta 2^{\alpha}>1$.
Let $(G,f)$ be $\alpha$-critical, where $G$ is house-free. Then
there is a subset $X\subseteq V(G)$
and a good function
$g$ on $G[X]$ with $g(v)\le f(v)$ for each $v\in X$,
such that $g^\alpha(X)\ge \delta f^\alpha(G)$ and $g^\alpha(N_G(v)\cap X)< \vare g^\alpha(X)$ for every vertex $v\in X$.
\end{thm}
\Proof By rational approximation, we may assume that
$f^\alpha$ is rational. Choose an integer $T>0$ such that $Tf^\alpha(v)$ is an integer for all $v\in V(G)$. Let $G'$ be obtained from $G$
by replacing each vertex $v$ by a clique $W_v$ of cardinality $Tf^\alpha(v)$, where
\begin{itemize}
\item the sets $W_v\;(v\in V(G))$ are pairwise disjoint;
\item for all distinct $u,v\in V(G)$ adjacent in $G$, $W_u$ is complete to $W_v$ in $G'$; and
\item for all distinct $u,v\in V(G)$ nonadjacent in $G$, $W_u$ is anticomplete to $W_v$ in $G'$.
\end{itemize}
It follows that $G'$ is also house-free, and $|G'|=Tf^\alpha(G)>T\ge 1$.
From \ref{nohouse} applied to $G'$, we deduce that either
\begin{itemize}
\item there are disjoint sets $A',B'\subseteq V(G')$, complete to each other, with $|A'|,|B'|\ge \vare\delta |G'|$, or
\item there exists $X'\subseteq V(G')$ with $|X'|\ge \delta|G'|$ such that $G'[X']$ has maximum degree less than $\vare|X'|$.
\end{itemize}
Suppose that the first bullet holds. Let $A$ be the set of vertices in $G$ such that $W_v\cap A'\ne \emptyset$, and define
$B$ similarly. Then $A$ is complete to $B$ in $G$. Moreover
$$f^\alpha(A)\ge |A'|/T\ge \vare\delta |G'|/T\ge \vare\delta f^\alpha(G)>\vare\delta,$$
and similarly $f^\alpha(B)\ge \vare\delta f^\alpha(G)$. By \ref{strongEH}, $\vare\delta \le 2^{-\alpha}$, contrary
to the hypothesis.

Thus the second bullet holds. Let $X$ be the set of all $v\in V(G)$ such that $W_v\cap X'\ne \emptyset$; and for each $v\in V(G)$
let $g(v)$
satisfy $T(g(v))^\alpha=|W_v\cap X'|$. Thus $g^\alpha(X)=|X'|\ge \delta|G'|=\delta f^\alpha(G)$, and $g(v)\le f(v)$
for each $v\in V(G)$. Let $v\in X$.
The union of the sets $W_u\cap X'$ over all $u\in N(v)\cap X$ has cardinality less than $\vare|X'|$ (indeed, less than
$\vare|X'|-|W_v\cap X'|+1$); and so $Tg^\alpha(N(v)\cap X)< \vare|X'|= \vare Tg^\alpha(X)$, that is,
$g^\alpha(N(v))< \vare g^\alpha(X)$. This proves \ref{fnohouse}.~\bbox

Since $P_4$-free graphs are perfect, a theorem of Erd\H{o}s and Hajnal~\cite{EH} (see also Alon, Pach and Solymosi~\cite{alon})
implies:
\begin{thm}\label{noforcer}
There exists $\vare>0$ such that if $G$ is forcer-free then $G$ has a clique or stable set of cardinality at least $|G|^\vare$.
\end{thm}

We also need a theorem of Jacob Fox (he did not publish his proof, but we gave a proof in~\cite{pathsandanti}):
\begin{thm}\label{Fox2}
Let $H$ be a graph for which there exists a constant $\delta>0$
such every $H$-free graph $G$ has a clique or stable set of cardinality at least  $|G|^{\delta}$.
Then
every $H$-free graph is $\frac{3}{\delta}$-narrow.
\end{thm}

By combining \ref{noforcer} and \ref{Fox2} we obtain:
\begin{thm}\label{forcernarrow}
There exists $\alpha\ge 1$ such that every forcer-free graph is $\alpha$-narrow.
\end{thm}

We deduce:
\begin{thm}\label{homog}
Let $\alpha'\ge 1$ such that every forcer-free graph is $\alpha'$-narrow.
Let $\alpha\ge \alpha'$, and let $G$ be a graph such that every proper induced subgraph is $\alpha$-narrow.
Let $g$ be a good function on $G$. Let $Z\subseteq V(G)$ with $|Z|>1$, such that
$G[Z]$ is connected and guarded.
Let $d$ be the maximum of $g^\alpha(N_G(v)\cap Z)$ over all $v\in Z$.
Then $g^\alpha(Z)\le \max(2d,d^{1-\alpha'/\alpha})$.
\end{thm}
\Proof
If $G[Z]$ is not anticonnected, there are two vertices $u,v\in Z$ such that $N(u)\cup N(v)=Z$, and so
$g^\alpha(Z)\le 2 d$ as required.
So we may assume that $G[Z]$ is anticonnected.
Let us list all subsets $X$ of $Z$ with the properties that $X$ is a homogeneous set of $G[Z]$, and $X\ne Z$, and $X$ is maximal
with these two properties; let these subsets be $W_1\ll W_k$ say. Thus $W_1\cup\cdots\cup W_k=Z$, because $|Z|\le 2$ and so
each singleton subset of $Z$ is a subset of one of $W_1\ll W_k$.

We claim that $W_1\ll W_k$ are pairwise disjoint. Suppose that $W_1\cap W_2\ne \emptyset$ say. Choose $w_1\in W_1\setminus W_2$
and $w_2\in W_2\setminus W_1$. If $w_1,w_2$ are nonadjacent, then since $W_2$ is a homogeneous set, $w_1$ has no neighbours in $W_2$,
and so, since $W_1$ is homogeneous, each vertex of $W_2$ has no neighbour in $W_1$; and so $G[W_1\cup W_2]$ is not connected. If $w_1,w_2$ are adjacent,
then similarly $G[W_1\cup W_2]$ is not anticonnected. Since $G[Z]$ is both connected and anticonnected, it follows that
$W_1\cup W_2\ne Z$. But $W_1\cup W_2$ is a homogeneous set of $G[Z]$, contrary to the maximality of $W_1$. This proves that
$W_1\ll W_k$ form a partition of $Z$, and so $k>1$.

Choose $w_i\in W_i$ for $1\le i\le k$, and let
$G'$ be the graph induced on $\{w_1\ll w_k\}$. From the hypothesis, $G'$ is forcer-free, and so $\alpha'$-narrow.
Let $t=\alpha/\alpha'$ and $r=2^{-1/\alpha'}$; thus $r^{\alpha'} = 1/2$, and $r^{\alpha}= 2^{-t}$.
For $1\le i\le k$, let $P_i$ be a perfect subgraph of $G[W_i]$ with $g(P_i)$ maximum, and let $p(w_i)=g(P_i)$.
By an application of Lov\'asz' substitution lemma~\cite{lovasz}, it follows that $p$ is a good function on $G'$, and so $p^{\alpha'}(G')\le 1$. For $1\le i\le k$, $G[W_i]$ is $\alpha$-narrow
(because $k>1$). Since $g(v)/p(w_i)\;(v\in W_i)$
is good on $G[W_i]$, it follows that $g^\alpha(W_i)\le p(w_i)^\alpha$. But also, since $G[Z]$ is connected, and
$W_i\ne Z$, there is a vertex $v\in Z\setminus W_i$ complete to $W_i$; and so $g^\alpha(W_i)\le d$ by hypothesis.
Thus
$$g^\alpha(W_i) \le \min(p(w_i)^\alpha,d)\le p(w_i)^{\alpha'}d^{1-1/t}.$$
Hence
$$g^\alpha(Z)=\sum_{1\le i\le k}g^\alpha(W_i)\le d^{1-1/t}\sum_{1\le i\le k}p(w_i)^{\alpha'}= d^{1-1/t}p^{\alpha'}(G')\le d^{1-1/t}.$$
This proves \ref{homog}.~\bbox

We deduce \ref{mainthm3}, which we restate:
\begin{thm}\label{mainthm4}
There exists $\alpha\ge 1$ such that for every $\alpha$-critical pair $(G,f)$, there is a hole-with-hat in $G$.
\end{thm}
\Proof
Let $\vare=1/6$, and choose $\delta>0$ satisfying \ref{nohouse}. From \ref{noforcer}, there exists $\alpha'\ge 1$ such that 
every forcer-free graph is 
$\alpha'$-narrow. Let $\alpha\ge 2$, such that  $\vare\delta 2^{\alpha/\alpha'}>1$.
We claim that $\alpha$ satisfies the theorem. Suppose not; then there is an 
$\alpha$-critical pair ($G,f)$, such that $G$ is hole-with-hat-free. 
By \ref{fnohouse}, there is a subset $X\subseteq V(G)$
and a good function
$g$ on $G[X]$ with $g(v)\le f(v)$ for each $v\in X$,
such that $g^\alpha(X)\ge \delta f^\alpha(G)>\delta$ and $g^\alpha(N_G(v)\cap X)< \vare g^\alpha(X)$ for every vertex $v\in X$.
Let $g^\alpha(X)=\lambda$; then $\delta\le \lambda\le f^\alpha(G)$. 
Let $H=G[X]$, and define $w(v)=g^\alpha(v)/\lambda$ for each $v\in X$. Then $(H, w)$ is a weighted graph. 
\\
\\
(1) {\em $(H,w)$ is $\vare$-coherent.}
\\
\\
By \ref{smalldeg}, for each $v\in V(H)$, 
$g(v)\le f(v)< 1-4^{-1/\alpha}\le  1/2$, and so
$w(v)< (1-4^{-1/\alpha})^\alpha/\lambda \le 2^{-\alpha}/\lambda\le \vare$ (since $\lambda\ge \delta$).
Also, for each vertex $v\in V(H)$, 
$$\lambda w(N_H(v))=g^\alpha(N_G(v)\cap X) < \vare g^\alpha(X)=\vare\lambda,$$
and so $w(N_H(v)) <\vare.$
Third, by \ref{strongEH}, if $A,B\subseteq V(H)$ are disjoint and anticomplete, then
not both $w(A), w(B)\ge \vare$, since $\vare>  2^{-\alpha}/\lambda$ (because $\lambda\ge \delta$). 
Consequently $(H, w)$ is $\vare$-coherent.
This proves~(1).

\bigskip

By \ref{flat}, and since $\vare=1/6$, there exists $Z\subseteq V(H)$ with $|Z|>1$, such that
$H[Z]$ is connected and guarded, and
$w(Z)>1-4\vare$. 
Let $d$ be the maximum of $\lambda w(N_G(v)\cap Z)$ over all $v\in Z$. Hence $d\le \vare\lambda$.
By \ref{homog}, 
$$\lambda w(Z)\le \max(2d,d^{1-\alpha'/\alpha})\le \max(2\vare\lambda,(\vare\lambda)^{1-\alpha'/\alpha}).$$ 
But $2\vare\lambda\ge (\vare\lambda)^{1-\alpha'/\alpha}$, since $2(\vare\lambda)^{\alpha'/\alpha}\ge 1$ (since $\lambda\ge \delta$). Thus
$\lambda w(Z)\le 2\vare\lambda$, and so $w(Z)\le 2\vare$, contradicting that $w(Z)> 1-4\vare$ and $\vare=1/6$.
This proves \ref{mainthm4}.~\bbox

\end{document}